\documentclass[onecolumn]{elsarticle}
\makeatletter
\def\ps@pprintTitle{%
	\let\@oddhead\@empty
	\let\@evenhead\@empty
	\def\@oddfoot{\centerline{\thepage}}%
	\let\@evenfoot\@oddfoot}
\makeatother

\usepackage{amsfonts}
\usepackage{amsmath}
\usepackage{amscd}
\usepackage[active]{srcltx} 
\usepackage{amssymb, color}
\usepackage[latin1]{inputenc}
\usepackage{enumerate}

\newtheorem{thm}{Theorem}[section]
 \newtheorem{cor}[thm]{Corollary}
 \newtheorem{lema}[thm]{Lemma}
 \newtheorem{prop}[thm]{Proposition}
 
 \newdefinition{Def}{Definition}
 
 \newenvironment{dem}{\vspace{1ex}\noindent{ Proof. }\hspace{0.5em}}
 {\hfill\qed\vspace{1ex}}
\newtheorem{example}{Example}
\numberwithin{equation}{section}

\newcommand{\HH}{\mathcal{H}}
\newcommand{\St}{\mathcal{S}}

\newcommand{\M}{\mathcal{M}}
\newcommand{\N}{\mathcal{N}}
\newcommand{\Q}{\mathcal{Q}}

\newcommand{\PI}[2]{\left\langle \,#1 , #2\, \right\rangle}

\newcommand{\mc}[1]{\mathcal{#1}}
\newcommand{\ol}{\overline}

\newcommand{\ra}{\rightarrow}

 \DeclareMathOperator{\tr}{tr}

\DeclareMathOperator{\rg}{R}

\begin{document}

\begin{frontmatter}

\title{Global solutions of approximation problems in Hilbert spaces}

\author[FI,IAM]{Maximiliano Contino \corref{ca}}
\ead{mcontino@fi.uba.ar}
\author[IAM]{Mar\'ia Eugenia Di Iorio y Lucero}
\ead{eugenia.diiorio@gmail.com}
\author[CIFASIS]{Guillermina Fongi}
\ead{fongi@cifasis-conicet.gov.ar}
\cortext[ca]{Corresponding author}

\address[FI]{%
Facultad de Ingenier\'{\i}a, Universidad de Buenos Aires\\
Paseo Col\'on 850 \\
(1063) Buenos Aires,
Argentina 
}

\address[IAM]{%
Instituto Argentino de Matem\'atica ``Alberto P. Calder\'on'' \\ CONICET\\
Saavedra 15, Piso 3\\
(1083) Buenos Aires, 
Argentina }

\address[CIFASIS]{%
	Centro Franco Argentino de Ciencias de la Informaci\'on y de Sistemas\\
	CONICET\\
	Ocampo y Esmeralda\\
	(2000) Rosario, Argentina}

\begin{abstract}
We study three well-known minimization problems in Hilbert spaces: the weighted least squares problem and the related problems of abstract splines and smoothing. In each case we analyze the solvability of the problem for every point of the Hilbert space in the corresponding data set, the existence of an operator that maps each data point to its solution in a linear and continuous way and the solvability of the associated operator problem in a fixed $p$-Schatten norm.
\end{abstract}


\end{frontmatter}

\section{Introduction}

In this work we focus our attention on the following well-known approximation and interpolation problems: the weighted least squares problem and the related problems of abstract splines  and smoothing.
Weighted least squares problems were studied in \cite{[CorMae05]} and some applications include Sard's approximation, least squares and curve fitting processes on a closed subspace \cite{Sard},  signal processing \cite{EldarA, EldarB} and sampling theory \cite{AnteCorach,Smale}.
Regarding the abstract spline and smoothing problems, Atteia \cite{[Att65]} obtained an abstract formulation which resumed most of the spline-type functions. This theory was developed by many authors, see for example Anselone and Laurent \cite{[AnsLau68]}, Shekhtman \cite{[She80]}, de Boor \cite{[Boo81]}, Izumino \cite{[Izu83]} and the surveys by Champion, Lenard and Mills \cite{[ChaLenMil96],[ChaLenMil00]}. The abstract spline and smoothing problems were generalized to bounded linear operators in \cite{Spline} and have been applied in many areas, such as approximation theory, numerical analysis and statistics, among others. See, for example \cite{[Sch46_A]}, \cite{[Sch46_B]} and \cite{Holladay}.

We study the conditions under which these minimization problems admit a bounded global solution, i.e., when it is possible to guarantee not only the existence of solutions for every point (of a Hilbert space) but also the existence of an operator that assigns to each point a solution in a linear and continuous way. We also study operator versions of these problems when the $p$-Schatten norms are considered and we relate  the existence of solutions of such problems to the existence of bounded global solutions. 
Let us fix some notations: $\HH, \mc{E}, \mc{F}$ are separable complex Hilbert spaces,  $L(\HH, \mc{F})$  is the set of bounded linear operators from $\HH$ to $\mc{F}$ and $L(\HH):=L(\HH,\HH).$ 

\subsection*{The weighted least squares problem}

Given $A\in L(\mathcal H, \mc{F})$, $W\in L(\mathcal F)$ positive (semidefinite) and $x \in \mathcal F$, a weighted least squares solution (or $W$-LSS) of the equation $Az = x$  is a vector $u\in \HH$ such that
\begin{equation}\tag{WLSP} \label{WLSP}
\|Au-x\|_{W} \leq \|Az-x\|_{W}, \quad \textrm{ for every }z \in \HH,
\end{equation}
where $\Vert x \Vert_W=\PI{Wx}{x}$ is the seminorm associated to $W$. 


Our main concern is to determine conditions for the existence of global solutions of \eqref{WLSP}. First, we analyze the existence of solutions of  \eqref{WLSP} for every $x\in \mathcal F$. Also, we study the  existence of a \emph{bounded global solution} of \eqref{WLSP} or a  \emph{$W$-inverse of $A$}, i.e., when there exists an operator $G \in L(\mc{F},\HH)$ such that, for each $x \in \mc{F}$, $G x$ is a $W$-LSS of $Az=x,$ see \cite{[RaoMit73],[CorFonMae13]}.
Finally, we study the following \textit{operator least squares problem}: given $A\in L(\mathcal H,\mc{F})$ (not necessarily with closed range) and $W\in L(\mathcal F)$ positive such that $W^{1/2}\in S_p$ for some $1\leq p <\infty,$ analyze the existence of
\begin{equation}\tag{OWLSP} \label{eqp}
\min_{X\in L(\mc{F}, \mathcal H)} \|AX-I\|_{p, W},
\end{equation}
where $\Vert Y \Vert_{p,W}=\Vert W^{1/2}Y \Vert_p, \mbox{ for } Y\in L(\mc{F}).$ Problem \eqref{eqp} was studied in \cite{Contino} when $A \in L(\HH)$ is a  closed range operator.

\subsection*{Spline and smoothing problems}
A similar approach as in  the case of the weighted least squares problem  will be done for the spline and the smoothing problems. 
Consider  $T \in L(\HH,\mc{E})$ and $V \in L(\HH,\mc{F})$. 

The \emph{classical spline problem}: given  $f_0 \in R(V),$ determine whether there exists
$$\min \Vert Th \Vert, \mbox{ subject to } Vh=f_0.$$  
The \textit{classical smoothing problem}: given $f_0\in \mathcal F$, analyze if there exists
$$
\underset{h\in \HH}{\min }(\|Th\|^2+\|Vh-f_0 \|^2).
$$

For each of the above two problems we study the existence of different global solutions as in the case of the weighted least squares problem.

In conclusion, for each of the three problems presented, we give conditions for the existence of solutions for every point in the corresponding data set, the existence of an operator that provides the solution at every data point in a continous way and the solution of the operator version of each of these problems considering the $p$-Schatten norms. 
We compare the different approaches and establish necessary and sufficient conditions for the existence of the described global solutions. 

The paper is organized as follows. In Section 2
we recall the notion of compatibility between a positive operator and a closed subspace of $\mathcal F$. Also, we collect certain properties of the Schatten class operators that will be used along the paper.
In Section 3, we study  global solutions of the weighted least squares problem \eqref{WLSP}, where the range of $A$ is not necessarily closed. 
It is proved that $A$ admits a $W$-inverse if and only if $R(A)+N(W)$ is closed and the pair $(W, R(A)+N(W))$ is compatible. This generalizes the fact that the classical least squares problem $\underset{z\in \mathcal {H}}{\min}{\|Az-x\|}$ admits a solution for every $x\in \mathcal F$ if and only if $R(A)$ is closed. Also, we show that \eqref{WLSP} admits a solution for every $x \in \mc{F}$ if and only if $A$ admits a $W$-inverse, or equivalently, if \eqref{eqp} admits a solution.

In Section 4, conditions for the existence of bounded global solutions in both the spline and the smoothing problems are given.
The solutions of these problems are compared to the solutions of the corresponding operator problems in the $p$-Schatten class. Finally, $W$-optimal inverses are compared to $W$-inverses and another connection between \eqref{WLSP} and the smoothing problem is given.

\section{Preliminaries}

Throughout $\HH, \mc{E}, \mc{F}$  are separable complex Hilbert spaces, $L(\HH, \mc{F})$ is the set of bounded linear operators from $\HH$ to $\mc{F},$  $L(\HH):=L(\HH,\HH),$ and $L(\HH)^{+}$ denotes the cone of semidefinite positive operators. The symbol $\leq$ stands for the order in $L(\HH)$  induced by $L(\HH)^+$, i.e., given $A,B \in L(\HH)$, $A \leq B$ if $B-A \in L(\HH)^+.$
For any $A \in L(\HH,\mc{F}),$ its range and nullspace are denoted by $R(A)$ and $N(A)$, respectively. Finally,  $A^{\dagger}$ denotes the  Moore-Penrose inverse of the operator $A \in L(\HH,\mc{F}).$ 

Given two closed subspaces $\M$ and $\N$ of $\HH,$ $\M \dot{+} \N$ denotes the direct sum of $\M$ and $\N$. Moreover, $\M \oplus \N$ stands for their (direct) orthogonal sum and $\M \ominus \N:= \M \cap (\M\cap \N)^{\perp}.$ 

Consider also, the classical inner product on $\mc{E} \oplus \mc{H}$ 
\begin{equation} \label{pi}
\PI{(e,h)}{(e',h')}=\PI{e}{e'} + \PI{h}{h'}, \mbox{ for } e,e' \in \mc{E}, h, h' \in \mc{H},
\end{equation}
together with the associated norm $\Vert (e,h) \Vert^2 = \PI{e}{e} + \PI{h}{h}, \mbox{ for } (e,h) \in \mc{E} \oplus \mc{H}.$

If $\HH$ is decomposed as a direct sum of closed subspaces $\HH=\M \dot{+} \N,$ the projection onto $\M$ with nullspace 
$\N$ is denoted by $P_{\M {\mathbin{\!/\mkern-3mu/\!}} \N}$ and $P_{\M} = P_{\M {\mathbin{\!/\mkern-3mu/\!}} \M^{\perp}}.$ Also, $\Q$ denotes the subset of $L(\HH)$ of oblique projections, i.e. $\Q=\{Q \in L(\HH): Q^{2}=Q\}.$ 

Given $W\in L(\mc{F})^{+}$ and a (non necessarily closed) subspace $\St$ of $\mc{F},$ the $W$-orthogonal complement of $\St$ is $\St^{\perp_{W}}=\{x \in \mc{F}: \PI{Wx}{y}=0, \ y\in \St\}=W^{-1}(\St^{\perp})=W(\St)^{\perp}.$
If $\St$ is a closed subspace of $\mc{F},$ the pair $(W,\St)$ is $\it{compatible}$ if there exists $Q\in \Q$ with $R(Q)=\St$ such that $WQ=Q^{*}W.$  

The next proposition, proved in \cite[Prop.~3.3]{CMSSzeged}, characterizes the compatibility of the pair $(W,\St).$

\begin{prop} \label{Comp 1} 
	Consider  $W\in L(\mc{F})^{+}$ and a closed subspace $\St \subseteq \mc{F}$. Then   the pair $(W,\St)$ is compatible if and only if $\mc{F}=\St+\St^{\perp_{W}}$.
\end{prop}

The notion of Schur complement of $W$ to $\St,$ for an operator $W\in L(\mc{F})^{+}$ and a closed subspace $\St \subseteq \mc{F},$ was introduced by M. G. Krein in \cite{Krein} and later rediscovered by Anderson and Trapp \cite{Shorted2}. They proved that the set
$\{ X \in L(\mc{F}): \ 0\leq X\leq W \mbox{ and } R(X)\subseteq \St^{\perp}\}$ has a maximum element. The Schur complement of $W$ to $\St$ is defined by
$$W_{/\St} := \max \ \{ X \in L(\mc{F}): \ 0\leq X\leq W \mbox{ and } R(X)\subseteq \St^{\perp}\}.$$

\vspace{0,5cm}
Let $T\in L(\HH, \mc{E})$ be a compact operator. By $\{\lambda_k(T)\}_{k\geq1}$ we denote the eigenvalues of $\vert T \vert := (T^{*}T)^{1/2}\in L(\HH),$ where each eigenvalue is repeated according to its multiplicity. Let $1\leq p < \infty,$ we say that $T$ belongs to the $p$-Schatten class $S_p(\HH,\mathcal E),$ 
if $\sum_{k\geq1}^{} \lambda_{k}(T) ^{p}<\infty$ and, the $p$-Schatten norm is given by 
$\Vert T\Vert_p:= (\sum_{k\geq1}^{} \lambda_{k}(T) ^{p})^{1/p}.$ For short, we write $S_p:=S_p(\HH,\HH).$ The set $S_p(\HH,\mc{E})$ is a vector space and $T \in S_p(\HH, \mc{E})$ if and only if $\vert T \vert \in S_p$ (see \cite[Theorem 7.6 and Theorem 7.8]{[Wei12]}).
If $T \in S_p(\HH,\mc{E}),$ then $$\Vert T \Vert_p^p= \tr ( \vert T \vert^p),$$ where $\tr$ denotes the trace of an operator.

Observe that, $\vert T \vert X \in S_p$ for every $X \in L(\HH)$ and $\Vert TX \Vert_p = \Vert \vert T \vert X \Vert_p.$ The reader is referred to \cite{Ringrose,Simon,[Wei12]} for further details on these topics.

\medskip

The following result will be useful along this paper, for its proof see \cite[Proposition 2.9]{Contino}. A more general result can be found in \cite[Proposition 2.5]{Nashed}.
\begin{prop}\label{Prop Nashed} Let $S, T \in S_p$ for some  $1\leq p < \infty.$ If $T^{*}T\leq S^{*}S \mbox{ then } \Vert T \Vert_p \leq \Vert S \Vert_p.$ \end{prop}

\vspace{0,5cm}
The Fr\'echet derivative will be instrumental to prove some results. We recall that, for a Banach space $(\mc{E}, \Vert \cdot \Vert)$ and an open set $\mathcal U \subseteq \mc{E},$ a function $F: \mc{E} \rightarrow \mathbb{R}$ is said to be \emph{Fr\'echet differentiable} at $X_0 \in \mathcal U$ if there exists $DF(X_0)$ a bounded linear functional such that
$$\lim\limits_{Y\rightarrow 0} \frac{|F(X_0+Y)-F(X_0) - DF(X_0)(Y)|}{\Vert Y \Vert}=0.$$ 
If $F$ is Fr\'echet differentiable at every $X_0 \in \mc{E}$, $F$ is called Fr\'echet differentiable  on $\mc{E}$  and the function $DF$ which assigns to every point $X_0 \in \mc{E}$ the derivative $DF(X_0),$ is called the Fr\'echet derivative of the function $F.$ If, in addition, the derivative $DF$ is continuous, $F$ is  said to be a \emph{class $\mc{C}^1$-function}, in symbols, $F \in \mc{C}^1(\mc{E}, \mathbb{R}).$

\begin{thm} \label{TeoD} Let $G_p: S_p \rightarrow \mathbb{R}^{+},$ $1 < p<\infty,$ $G_p(X)=\Vert X \Vert_p^{p},$ and let $X, Y \in S_p$.  Then, for $1< p<\infty,$ $G_p$ has a Fr\'echet derivative given by 
	$$DG_p(X)(Y)=p \ Re \ [ \tr ( \vert X \vert^{p-1} U^{*}Y)],$$ 
	where $Re(z)$ is the real part of a complex number $z$ and $X=U\vert X\vert$ is the polar decomposition of the operator $X,$ with $U$ the partial isometry such that $N(U)=N(X).$\end{thm}

\begin{dem} See \cite[Theorem 2.1]{Aiken}.
\end{dem}

Observe that, if $p=2.$ Then, $G_2(X)=\tr(X^*X)$ and $$DG_2(X)(Y)=2 \ Re \ \tr ( X^{*}Y), \ X, Y \in S_2.$$
\section{Global solutions of weighted least squares problems}\label{section3}

Given $A\in L(\mathcal H, \mc{F})$, $W\in L(\mathcal F)^+$ and $x \in \mathcal F$, a weighted least squares solution (or $W$-LSS) of the equation $Az = x$  is a vector $u\in \HH$ such that
\begin{equation}\tag{WLSP} 
\|Au-x\|_{W} \leq \|Az-x\|_{W}, \quad \textrm{ for every }z \in \HH,
\end{equation}
where $\Vert x \Vert_W=\PI{Wx}{x}$ is the seminorm associated to $W$. The related problem is the \textit{classical weighted least squares problem}.

If $W=I$, then problem \eqref{WLSP} is the well-known least squares problem. Given $A\in L(\HH, \mc{F})$ and  $x \in \mc{F},$ it can be proved  that  $u$ is a least squares solution of $Az =x$ if and only if $Au - x \in R(A)^{\perp}$ (see \cite[Theorem 3.1]{NashedLSS}).  As a consequence, it is not difficult to prove the following result.

\begin{prop} \label{propWLSS} Let $A\in L(\HH, \mc{F}),$ $W\in L(\mathcal F)^+$ and $x \in \mathcal F.$ Then $u$ is a $W$-LSS of $Az =x$ if and only if $Au - x \in W(R(A))^{\perp}$ or, equivalently, $A^*W(Au-x)=0.$ 
\end{prop}

\smallskip

To study the existence of solutions  of problem \eqref{WLSP} for every $x\in \mathcal F$ in the finite dimensional case, Rao and Mitra  introduced the notion of $W$-inverse \cite{[RaoMit73]}. Later on, the $W$-inverse was studied for operators in \cite{[CorFonMae13]} and \cite{Contino}.

\begin{Def}\rm Given $A\in L(\mathcal H, \mc{F})$ and $W\in L(\mathcal F)^+.$ An operator $G \in L(\mc{F},\HH)$ is called a \emph{$W$-inverse} of $A$ (or a \emph{bounded global solution} of problem \eqref{WLSP}) if for each $x \in \mc{F}$, $G x$ is a $W$-LSS of $Az=x$, i.e.
	$$
	\| AGx-x\|_W \leq \|Az-x\|_W, \quad \textrm{ for every }z \in \HH.
	$$
\end{Def}

Observe that $A$ has a $W$-inverse if the problem \eqref{WLSP} admits a solution for every $x\in \mathcal F$ and, moreover, it is possible to assign a $W$-LSS to each $x \in \mc{F}$ in a linear and continuous way.
\medskip

The following result gives necessary and sufficient conditions for problem \eqref{WLSP} to admit a bounded global solution, when $A$ is not  necessarily a closed range operator (cf. \cite{Contino}).

\begin{thm}\label{sol global sii  A-inversa} Let $A\in L(\mathcal H, \mc{F})$ and $W\in L(\mathcal F)^+.$ Then the following statements are equivalent:
	\begin{enumerate}
		\item[i)] $Az=x$ admits a $W$-LSS for every $x \in \mc{F},$ 
		\item[ii)]  $R(A)+W(R(A))^\perp=\mathcal{F},$
		\item[iii)] the normal equation \begin{equation} \label{ecuacion normal}
		A^*W(AX-I)=0
		\end{equation} admits a solution,
		\item[iv)] the operator $A$ admits a $W$-inverse.
	\end{enumerate}
	In this case, the set of $W$-inverses of $A$ is the set of solutions of \eqref{ecuacion normal}.
\end{thm}

\begin{dem} Along the proof we will use that $W(R(A))^{\perp}=N(A^*W).$
	
	$i) \Rightarrow ii):$ By Proposition \ref{propWLSS}, Problem \eqref{WLSP} admits a solution for every $x\in \mc{F}$ if and only if there exists $u$ such that $A^*WAu=A^*Wx$ for every $x \in \mathcal F$ or, equivalently $R(A^*W)=R(A^*WA)$. Then
	$$\mc{F}=(A^*W)^{-1}(R(A^*W))=(A^*W)^{-1}(A^*W(R(A)))=R(A)+W(R(A))^\perp.$$
	
	$ii) \Rightarrow iii):$ Suppose that $R(A)+W(R(A))^\perp=\mathcal{F}$. Then 
	$R(A^*W)=R(A^*WA)$ and the assertion follows by applying Douglas' Lemma \cite{Douglas}.
	
	$iii) \Rightarrow iv):$ There exists $X_0\in L(\mc{F},\HH)$ such that $A^*WAX_0=A^*W$ if and only if  $A^*W(AX_0x-x)=0$ for every $x\in \mathcal F$ or $AX_0x-x\in W(R(A))^\perp$ for every $x\in \mathcal F.$  Therefore, by Proposition \ref{propWLSS}, $X_0$ is a $W$-inverse of $A$.
	
	$iv) \Rightarrow i):$ It is straightforward.
	
	In this case, we have also proved that the set of $W$-inverses of $A$ is the set of solutions of \eqref{ecuacion normal}.
\end{dem}

\smallskip
\begin{cor} Let $A\in L(\mathcal H, \mc{F})$ and $W\in L(\mathcal F)^+.$  If $A$ admits a $W$-inverse then $(W, \overline{R(A)})$ is compatible.
\end{cor}
\begin{dem}
	If $A$ admits a $W$-inverse then, by Theorem \ref{sol global sii  A-inversa} and the identity $W(R(A))^{\perp}=W(\ol{R(A)})^{\perp},$ we get that $$\mathcal{F}=R(A)+W(R(A))^\perp \subseteq \ol{R(A)} + W(\ol{R(A)})^\perp.$$  Therefore, by Proposition \ref{Comp 1}, the pair $(W, \overline{R(A)})$ is compatible.
\end{dem} 

A non closed range operator can admit a $W$-inverse, as the following example shows.
\begin{example}
	Consider $W\in L(\mathcal F)^+ $ with infinite dimensional nullspace and $A=P_{N(W)^\perp}+A_2P_{N(W)},$ where $A_2\in L(N(W))$ and $R(A_2)$ is not closed. Then $R(A)$ is not closed and 
	$$
	R(A)+W(R(A))^\perp=R(A)+W(N(W)^\perp)^\perp=R(A)+R(W)^\perp \supseteq N(W)^{\perp}+R(W)^\perp=\mathcal F,
	$$
	so that, by Theorem \ref{sol global sii  A-inversa}, it holds that $A$ admits a $W$-inverse. 
\end{example}

If $W=I$, it is well known that the least squares problem $Az=x$ admits a solution for every $x \in \mc{F}$
if and only if $R(A)$ is closed.
More generally:
\begin{prop} Let $A\in L(\mathcal H, \mc{F})$ and $W\in L(\mathcal F)^+.$ Then $A$ admits a $W$-inverse if and only if $R(A)+N(W)$ is closed and the pair $(W, R(A)+N(W))$ is compatible.
\end{prop}
\begin{dem} Suppose that $A$ admits a $W$-inverse. Then, by Theorem \ref{sol global sii  A-inversa}, it holds that $ \mc{F}=R(A)+W(R(A))^\perp=R(A)+W^{-1}((R(A))^\perp)$. Applying $W^{1/2}$, it follows that $R(W^{1/2})=W^{1/2}(R(A))+W^{1/2}(W^{-1}(R(A)^\perp))=W^{1/2}(R(A))+W^{1/2}(R(A))^\perp \cap R(W^{1/2})$.
	Therefore, $ W^{1/2}(R(A))$ is closed in $R(W^{1/2})$. Hence, $R(A)+N(W)=R(A)+N(W^{1/2})=W^{-1/2}(W^{1/2}(R(A))$ is closed.
	Since $\mc{F}=R(A)+W(R(A))^\perp=R(A)+N(W)+W(R(A)+N(W))^\perp,$ by Proposition \ref{Comp 1},  the pair $(W, R(A)+N(W))$ is compatible.
	Conversely, suppose that $R(A)+N(W)$ is closed and the pair $(W, R(A)+N(W))$ is compatible. Then $\mc{F}=R(A)+N(W)+W(R(A)+N(W))^\perp=R(A)+W(R(A))^\perp,$ because $N(W)\subseteq W(R(A))^\perp$. Therefore, by Theorem \ref{sol global sii  A-inversa}, $A$ admits a $W$-inverse.
\end{dem}

\begin{prop} \label{A inyectivo implica R(B) cerrado} Let $A\in L(\mathcal H, \mc{F})$ and $W\in L(\mathcal F)^+.$ If $A$ admits a $W$-inverse then $R(A)$ is closed if and only if  $R(A)\cap N(W)$ is closed.
\end{prop}
\begin{dem} If $A$ admits a $W$-inverse then, by Theorem \ref{sol global sii  A-inversa},
	it holds that $\mathcal F=R(A)+W(R(A))^\perp$.
	Since $W\in L(\mathcal F)^+$, it is not difficult to see that
	$R(A) \cap W(R(A))^\perp=R(A) \cap N(W).$ Suppose that
	$R(A)\cap N(W)$ is closed, then
	$\mathcal F=R(A)+W(R(A))^\perp=R(A)\dot{+}[W(R(A))^\perp\ominus(R(A)\cap N(W))].$
	Hence,  by \cite[Theorem 2.3]{[FilWil71]}, it follows that $R(A)$ is closed. The converse is straightforward.
\end{dem}

\subsection*{Operator weighted least squares problems}
In this subsection we are interested in studying weighted least squares problems for operators considering Schatten $p$-norms.

Given $A\in L(\mathcal H, \mc{F})$ and $W\in L(\mathcal F)^+$ such that $W^{1/2}\in S_p$ for some $1\leq p <\infty,$ 
the problem is to determine if there exists
\begin{equation}\tag{OWLSP} 
\min_{X\in L(\mc{F},\mathcal H)} \|AX-I\|_{p, W},
\end{equation}
where $\Vert Y \Vert_{p,W}=\Vert W^{1/2}Y \Vert_p, \mbox{ for } Y\in L(\mc{F}).$ 

We will refer to problem \eqref{eqp} as the $\emph{operator weighted least squares problem}.$

In order to study problem \eqref{eqp}, we introduce the following associated problem: given $A\in L(\mathcal H, \mc{F})$ and $W\in L(\mathcal F)^+,$ analyze the existence of 
\begin{equation} \label{eqop}
\min_{X\in L(\mc{F},\mathcal H)} (AX-I)^*W(AX-I),
\end{equation}
in the order induced in $L(\mc{F})$ by the cone of positive operators.
By studying problems \eqref{eqp} and \eqref{eqop} we will relate the existence of solutions of \eqref{eqp} to the existence of bounded global solutions of \eqref{WLSP}.

In \cite{Contino}, problems  \eqref{eqp} and \eqref{eqop} were studied for $A\in L(\mathcal H)$  such that $R(A)$ is closed. The results obtained in \cite{Contino} are also valid in the general case.

\begin{prop} \label{PropA} Let $A\in L(\mathcal H, \mc{F})$ and $W\in L(\mathcal F)^+$ such that $W^{1/2}\in S_p$ for some $1\leq p <\infty$. Then the following statements are equivalent:
	\begin{enumerate}
		\item [i)]  there exists  $\underset{X\in L(\mc{F},\mathcal H)}{\min} \|AX-I\|_{p, W},$
		\item [ii)] $R(A)+W(R(A))^\perp=\mathcal{F},$
		\item [iii)] there exists $\underset{X\in L(\mc{F},\mathcal H)}{\min} (AX-I)^*W(AX-I).$ 
	\end{enumerate}
	
	In this case, $$\min_{X\in L(\mc{F},\mathcal H)} \|AX-I\|_{p, W}= \Vert W_{ / \ol{R(A)}}^{1/2} \Vert_p \mbox { and }\min_{X\in L(\mc{F},\mathcal H)} (AX-I)^*W(AX-I)=W_{ / \ol{R(A)}}.$$   
\end{prop}

\begin{dem} The equivalence between $i),$ $ii)$ and $iii)$ follows by similar arguments as those in the proofs of \cite[Theorem 4.3]{Contino}  and \cite[Theorem 4.5]{Contino}, using Proposition \ref{propWLSS}.
	
	In this case. Let $X_0$ be a solution of problem \eqref{eqop}, i.e., $(AX_0-I)^*W(AX_0-I)=\underset{X\in L(\mc{F},\mathcal H)}{\min} (AX-I)^*W(AX-I).$  Then, in particular, $0\leq (AX_0-I)^*W(AX_0-I) \leq W$ and, by similar arguments as those found in \cite[Proposition 4.4]{Contino}, $A^*W(AX_0-I)=0.$ Therefore, since $A^*[(AX_0-I)^*W(AX_0-I)]=A^*X_0^*A^*W(AX_0-I)-A^*W(AX_0-I)=0,$ we have that $R((AX_0-I)^*W(AX_0-I))\subseteq R(A)^{\perp}.$ 
	Let $Z \in L(\mc{F})^+$ such that $Z \leq W$ and $R(Z) \subseteq R(A)^{\perp}.$Then, $$Z=(AX_0-I)^*Z(AX_0-I) \leq (AX_0-I)^*W(AX_0-I).$$ Therefore  $\underset{X\in L(\mc{F},\mathcal H)}{\min} (AX-I)^*W(AX-I)=(AX_0-I)^*W(AX_0-I)=\max \  \{ Z \in L(\mc{F}): \ 0\leq Z \leq W \mbox{ and } R(Z)\subseteq R(A)^{\perp}\} = W_{ / \ol{R(A)}}.$
	
	Finally, by Proposition \ref{Prop Nashed}, $\Vert AX_0-I \Vert_{p,W}= \Vert W_{ / \ol{R(A)}}^{1/2} \Vert_p$ and
	$$\min_{X\in L(\mc{F},\mathcal H)} \|AX-I\|_{p, W}=\Vert AX_0-I \Vert_{p,W}= \Vert W_{ / \ol{R(A)}}^{1/2} \Vert_p.$$ \end{dem}

\medskip
The next corollary summarizes the results of the section.

\begin{cor} \label{CorA} Let $A\in L(\mathcal H, \mc{F})$ and $W\in L(\mathcal F)^+$ such that $W^{1/2}\in S_p$ for some $1\leq p <\infty$. Then the following statements are equivalent:
	\begin{enumerate}
		\item [i)] $Az=x$ admits a $W$-LSS for every $x \in \mc{F},$ 
		\item [ii)] $A$ admits a $W$-inverse, i.e., for every $x \in \mc{F},$ $Az=x$ admits a $W$-LSS $Gx$ with $G \in L(\mc{F},\HH),$
		\item [iii)] there exists  $\underset{X\in L(\mc{F},\mathcal H)}{\min}  \|AX-I\|_{p, W}.$
	\end{enumerate}
\end{cor}

\smallskip

\section{Global solutions of spline and  smoothing problems}
Following similar ideas as those presented in Section \ref{section3}, we will study under which conditions the classical spline and smoothing problems admit global solutions. Moreover, we will relate bounded global solutions to the solutions of the associated operator minimization problems.

\subsection*{Splines problems}
Let $T \in L(\HH,\mc{E}),$ $V \in L(\HH,\mc{F})$ and $f_0 \in R(V)$, we study the existence of
\begin{equation} \label{spline11}
\min \Vert Th \Vert, \mbox{ subject to } Vh=f_0. 
\end{equation}

Suppose that $Vh_0=f_0,$ problem \eqref{spline11} is equivalent to study when the set
\begin{equation} \tag{SP}\label{spline}
\begin{split} 
sp(T,N(V),h_0)= \{ h \in h_0+ N(V) : \Vert Th \Vert=\underset{z \in N(V)}{\min} \Vert T(h_0+z)\Vert\}
\end{split}
\end{equation}
is not empty. 
We will refer to problem \eqref{spline} as the \emph{classical spline problem} and any element of the set $sp(T,N(V),h_0)$ is an \emph{abstract spline} or a $(T,N(V))$-spline interpolant to $h_0.$  

In order to obtain solutions of \eqref{spline} that depend continuously on $h$ we give the following definition.

\begin{Def} Let $T \in L(\HH,\mc{E})$ and $V \in L(\HH,\mc{F}).$ An operator $G \in L(\HH)$ is a \emph{bounded global solution} of \eqref{spline}  if
	\begin{equation} \label{splineglobal}
	Gh \in sp(T,N(V),h) \mbox{ for every } h \in \HH.
	\end{equation}
\end{Def}
\smallskip

We are also interested in comparing the bounded global solution of \eqref{spline} to the \emph{operator spline problem}: given $T \in S_p(\HH, \mc{E})$ for some $1 \leq p < \infty,$ $V \in L(\HH,\mc{F})$ and $B_0 \in L(\HH,\mc{F})$ such that $R(V) \subseteq R(B_0),$ 
analyze the existence of
\begin{equation}\tag{OSP}\label{uno}
\min_{VX=B_0} \|TX\|_{p},
\end{equation}
where $X \in L(\HH).$ 

\vspace{0,3cm}
We begin by studying problem \eqref{uno}. The next result characterizes the existence of solutions of \eqref{uno}  and  describes  the operators where the minimum is attained.

\begin{prop} \label{teo1} Let $T \in S_p(\HH,\mc{E})$ for some $1 \leq p < \infty,$ $V \in L(\HH,\mc{F})$ and $B_0 \in L(\HH,\mc{F})$ such that $R(B_0) \subseteq R(V).$ Then the following statements are equivalent:
	\begin{itemize}
		\item [i)] there exists $\underset{VX=B_0}{\min} \|TX\|_{p},$ 
		\item [ii)] $R(V^{\dagger} B_0) \subseteq N(V) + \left[ T^{\ast}T(N(V))\right]^{\perp},$
		\item [iii)] the normal equation 	\begin{equation} \label{Normal2}
		P_{N(V)} T^{\ast}T (P_{N(V)}X+V^{\dagger} B_0)=0
		\end{equation} admits a solution.
	\end{itemize}	
	In this case, $$\min_{VX=B_0} \|TX\|_{p}=\| [(T^{\ast} T)_{/N(V)}]^{1/2} V^{\dagger}B_0\|_{p}=\Vert TX_0 \Vert_p,$$ where
	$X_0$ is any solution of equation \eqref{Normal2}.
\end{prop}

\begin{dem} $i) \Leftrightarrow ii):$ Note that if $VX=B_0$ then $V^{\dagger}VX=P_{N(V)^{\perp}}X=V^{\dagger}B_0 \in L(\HH)$ (see Douglas' Lemma \cite{Douglas}). Then, $X=P_{N(V)}X +V^{\dagger}B_0$ and, since $\Vert T X \Vert_p=\Vert X \Vert_{p,T^*T},$ we get that
	\begin{equation}
	\min_{VX=B_0} \|TX\|_{p} =	\min_{X \in L(\HH)}
	\| P_{N(V)}X + V^{\dagger} B_0 \|_{p,T^{\ast}T}.
	\end{equation}
	
	Then, by \cite[Theorem 4.5]{Contino}, problem \eqref{uno} admits a solution
	if only if $$\rg(V^{\dagger} B_0) \subseteq N(V) + \left[ T^{\ast}T(N(V))\right]^{\perp}.$$
	
	$ii) \Leftrightarrow iii):$ It follows from  \cite[Theorem 2.4]{Contino}.
	
	In this case, by \cite[Theorem 4.5]{Contino} and \cite[Theorem 2.4]{Contino}, $$\min_{VX=B_0} \|TX\|_{p}=\| [(T^{\ast} T)_{/N(V)}]^{1/2} V^{\dagger}B_0\|_{p}=\Vert TX_0 \Vert_p,$$ where
	$X_0$ is any solution of equation \eqref{Normal2}.
	
\end{dem}

\begin{prop} \label{prop11} Let $T \in S_p(\HH,\mc{E})$ for some $1 \leq p < \infty,$ $V \in L(\HH,\mc{F})$ and $B_0 \in L(\HH,\mc{F})$ such that $R(B_0) \subseteq R(V).$ 
	Then $X_0 \in L(\HH)$ is a solution of  \eqref{uno} if and only if $$X_0x \in sp(T,N(V), V^{\dagger}B_0x) \mbox{ for every } x \in \HH.$$
\end{prop}

\begin{dem} Suppose $X_0 \in L(\HH)$ is a solution of \eqref{uno}. Then there exists $Y_0 \in L(\HH)$ such that $X_0=P_{N(V)}Y_0+V^{\dagger}B_0$ and $X_0$ is a solution of \eqref{Normal2}. Then, $Y_0$ is a solution of \eqref{Normal2} too. So, by \cite[Proposition 4.4]{Contino}, 
	$$(P_{N(V)}Y_0+V^{\dagger} B_0)^*T^*T(P_{N(V)}Y_0+V^{\dagger} B_0) \leq (P_{N(V)}Y+V^{\dagger} B_0)^* T^*T (P_{N(V)}Y+V^{\dagger} B_0), $$ for every $Y \in L(\HH).$ Or, equivalently, 
	$$\Vert T(P_{N(V)}Y_0+V^{\dagger} B_0) x  \Vert \leq \Vert T(P_{N(V)}Y+V^{\dagger} B_0) x  \Vert, \mbox{ for every } x \in \HH \mbox{ and } Y \in L(\HH).$$
	Let $ z \in \HH$ be arbitrary. For every $x \in \HH \setminus \{ 0\},$ there exists $Y\in L(\HH)$ such that $z=Yx.$ Therefore
	$$\Vert T(P_{N(V)}Y_0+V^{\dagger} B_0) x  \Vert \leq \Vert T(P_{N(V)}z+V^{\dagger} B_0x)  \Vert, \mbox{ for every } x, z  \in \HH.$$ 
	Then $X_0x \in  V^{\dagger} B_0x+ N(V),$  
	$$\Vert T X_0 x  \Vert \leq \Vert T h  \Vert, \mbox{ for every } h  \in V^{\dagger} B_0x+ N(V)$$ and $X_0x \in sp(T,N(V), V^{\dagger}B_0x).$
	
	Conversely, suppose that $X_0x \in sp(T,N(V), V^{\dagger}B_0x) \mbox{ for every } x \in \HH.$ Then, for every $x \in \HH,$ $X_0x \in V^{\dagger} B_0x+ N(V)$ and 
	$$\Vert T X_0 x  \Vert \leq \Vert T h  \Vert, \mbox{ for every } h  \in V^{\dagger} B_0x+ N(V).$$ It follows that $VX_0=B_0$ and  
	$$\Vert T X_0 x  \Vert \leq \Vert T(P_{N(V)}z+V^{\dagger} B_0x)  \Vert, \mbox{ for every } x, z \in \HH.$$ In particular, given $Y \in L(\HH),$ consider $z=Yx.$ Then 
	$$\Vert T X_0 x  \Vert \leq \Vert T(P_{N(V)}Y+V^{\dagger} B_0) x \Vert  \mbox{ for every } x \in \HH \mbox{ and } Y \in L(\HH),$$ or, equivalently, $$X_0^*T^*TX_0 \leq (P_{N(V)}Y+V^{\dagger} B_0)^* T^*T (P_{N(V)}Y+V^{\dagger} B_0), \mbox{ for every } Y \in L(\HH).$$ Then, by Proposition \ref{Prop Nashed}, 
	$$\Vert TX_0 \Vert_p = \Vert \vert T \vert X_0 \Vert_p \leq
	\| \vert T \vert (P_{N(V)}Y + V^{\dagger} B_0) \|_{p}=	\| P_{N(V)}X + V^{\dagger} B_0 \|_{p,T^{\ast}T} 
	$$
	for every $ Y \in L(\HH).$
	Therefore, $X_0$ is a solution of \eqref{uno}.
\end{dem}

The following result gives necessary and sufficient conditions for the operator spline problem \eqref{uno} to have a solution for every $B_0 \in L(\HH,\mc{F}).$ Moreover, it shows that this is equivalent to the condition that guarantees the existence of a bounded global solution of the classical spline problem \eqref{spline}.

\begin{thm} Let $T \in S_p(\HH,\mc{E})$ for some $1 \leq p < \infty$ and $V \in L(\HH,\mc{F}).$ Then the following statements are equivalent:
	\begin{enumerate}
		\item [i)] there exists $\underset{VX=B_0}{\min} \|TX\|_{p}$ for every $B_0 \in L(\HH,\mc{F})$ such that $R(B_0) \subseteq R(V),$  
		\item [ii)] there exists a bounded global solution of \eqref{spline},
		\item [iii)] the pair $(T^{\ast}T, N(V))$ is compatible,
		\item [iv)]  $sp(T,N(V),h_0)$ is nonempty for every $h_0\in \HH$. 
	\end{enumerate} 
\end{thm}

\begin{dem} $i) \Leftrightarrow ii):$  Suppose that $X_0 \in L(\HH)$ is a solution of \eqref{uno} for $B_0=V.$ Consider $G:=X_0P_{N(V)^{\perp}} \in L(\HH).$ Then, by Proposition \ref{prop11}, 
	$$
	Gh=X_0(P_{N(V)^{\perp}}h) \in sp(T,N(V),P_{N(V)^{\perp}}h) \mbox { for every }  h \in \HH.
	$$
	Note that $sp(T,N(V),P_{N(V)^{\perp}}h)=sp(T,N(V),h)$, because $P_{N(V)^{\perp}}h+N(V)=h+N(V).$
	Hence,
	$$
	Gh \in sp(T,N(V),h) \mbox { for every }  h \in \HH,
	$$
	so that $G$ is a bounded global solution of \eqref{spline}.
	
	Conversely, suppose that $G \in L(\HH)$ is a bounded global solution of \eqref{spline} and $B_0 \in L(\HH,\mc{F}).$ Set $X_0:=GV^{\dagger}B_0 \in L(\HH),$ then
	$$X_0x=G(V^{\dagger}B_0x) \in sp(T,N(V),V^{\dagger}B_0x) \mbox{ for every } x \in \HH.$$ Therefore, by Proposition \ref{prop11}, $X_0$ is a solution of \eqref{uno}.
	
	$i) \Leftrightarrow iii):$ Suppose that \eqref{uno} has a solution for every $B_0 \in L(\HH,\mc{F}).$ Then, by Proposition \ref{teo1}, $$\rg(V^{\dagger} B_0) \subseteq N(V) + \left[ T^{\ast}T(N(V)\right]^{\perp}$$ for every $B_0$ such that $R(B_0) \subseteq R(V).$  
	Consider $B_0$ such that $R(B_0)=R(V),$ then 
	$N(V)^{\perp}= \rg(V^{\dagger} B_0) \subseteq N(V) + \left[ T^{\ast}T(N(V)\right]^{\perp},$ so that
	$\HH= N(V) + \left[ T^{\ast}T(N(V)\right]^{\perp}$ and, by Proposition \ref{Comp 1}, the pair $(T^{\ast}T, N(V))$ is compatible. 
	
	Conversely, let the pair $(T^{\ast}T, N(V))$ be compatible. By Proposition \ref{Comp 1}, $\HH= N(V) + \left[ T^{\ast}T(N(V)\right]^{\perp}.$ Then, for every $B_0,$ $\rg(V^{\dagger} B_0)  \subseteq N(V) + \left[ T^{\ast}T(N(V)\right]^{\perp}$ and, by Proposition \ref{teo1},  \eqref{uno} has a solution for every $B_0 \in L(\HH,\mc{F})$ such that $R(B_0) \subseteq R(V).$ 
	
	$iii) \Leftrightarrow iv):$ See \cite[Theorem 3.2]{Spline}.
\end{dem}

\subsection*{Smoothing problems}

Let $T \in L(\HH,\mc{E}),$ $V \in L(\HH,\mc{F})$ and $f_0\in \mathcal F.$  A problem that is naturally associated with \eqref{spline} is to find
\begin{equation} \tag{SMP} \label{CSP}
\underset{h\in \HH}{\min }(\|Th\|^2+\|Vh-f_0 \|^2).
\end{equation} 
We will refer to \eqref{CSP} as the \textit{classical smoothing problem} and its solutions are called \textit{smoothing splines}. 

As before, we also study the problem of finding a \emph{bounded global solution} of problem \eqref{CSP}; i.e., we analyze if there exists an operator $G\in L(\mathcal F, \HH)$ such that
\begin{equation} \label{GlobalCSP}
\|TGf\|^2+\|VGf-f\|^2=\underset{h\in \HH}{\min }(\|Th\|^2+\|Vh-f\|^2), \textrm{ for every } f \in \mathcal F.
\end{equation}

Several properties of bounded global solutions of problem \eqref{CSP} were given in \cite[Section 4]{[CorFonMae16]} for $V \in L(\mc{H},\mc{F})$ with closed range. 

We are interested in characterizing bounded global solutions of \eqref{CSP} in the general case and comparing them with the solutions of the  following \emph{operator smoothing} problem: given $T \in S_2(\HH,\mc{E}),$  $V \in S_2(\HH, \mc{F})$ and $B_0 \in S_2(\HH, \mc{F}),$  analyze the existence of
\begin{equation} \tag{OSMP}\label{dos}
\underset{X \in L(\HH)}{\min} \ (\Vert TX \Vert_{2}^2 + \Vert VX-B_0 \Vert_2^2).
\end{equation}

\smallskip

Define $K, B_0': \HH \ra \mc{E} \oplus \mc{F},$ 
\begin{eqnarray}
Kh=(Th,Vh) \mbox{ for } h \in \HH, \label{eqK} \\ 
B_0'h=(0,B_0h) \mbox{ for } h \in \HH \label{eqB}.
\end{eqnarray}
We will consider the inner product and the associated norm on $\mc{E} \oplus \mc{F}$ as in \eqref{pi}. It is straightforward to check that the adjoint of $K,$ $K^*: \mc{E} \oplus \mc{F}\ra \HH,$ is $K^*(e,f)=T^*h+V^*f, \mbox{ for } e \in \mc{E} \mbox{ and } f \in \mc{F}$ and the adjoint of $B_0',$ $B_0'^*: \mc{E} \oplus \mc{F}\ra \HH,$ is $B_0'^*(e,f)=B_0^*f, \mbox{ for } e \in \mc{E} \mbox{ and } f \in \mc{F}.$

\begin{lema}  \label{teo2} Let $T \in L(\HH,\mc{E}),$  $V \in L(\HH, \mc{F})$ and $B_0 \in L(\HH,\mc{F}).$ Set $K$ and $B_0'$ as in \eqref{eqK} and \eqref{eqB}. Then there exists $X_0\in L(\HH)$ such that 
	$$(KX_0-B_0')^*(KX_0-B_0')=\underset{X \in L(\HH)}{\min}(KX-B_0')^*(KX-B_0'),$$  where the order is the one induced in $L(\HH)$ by the cone of positive operators, if and only if $X_0$ is a solution of the normal equation
	\begin{equation} \label{Normal1}
	(T^*T+V^*V)X=V^*B_0.
	\end{equation}
	
\end{lema}

\begin{dem} This follows in a similar way as in the proof of Proposition \ref{PropA} and Theorem \ref{sol global sii  A-inversa} using the fact that $u$ is a least squares solution of the equation $Kz=B_0'x$ for every $x \in \HH$ if and only if $u$ is a solution of $K^{*}(Ky-B_0'x)=0,$ see Proposition \ref{propWLSS}.
\end{dem}

\begin{prop} \label{teo3} Let $T \in S_2(\HH,\mc{E}),$  $V \in S_2(\HH, \mc{F})$ and $B_0 \in S_2(\HH, \mc{F}).$ Then
	the following statements are equivalent:
	\begin{enumerate}
		\item [i)] there exists $\underset{X \in L(\HH)}{\min} \ (\Vert TX \Vert_{2}^2 + \Vert VX-B_0 \Vert_2^2),$
		\item [ii)] the normal equation $(T^*T+V^*V)X=V^*B_0$ admits a solution.
	\end{enumerate}
\end{prop}

\begin{dem} Let $K$ and $B_0'$ be as in \eqref{eqK} and \eqref{eqB}. If $X_0$ is a solution of the normal equation \eqref{Normal1}, then by Lemma \ref{teo2}, $(KX_0-B_0')^*(KX_0-B_0')\leq (KX-B_0')^*(KX-B_0'), \mbox{ for every }  X\in L(\HH).$ 	
	Observe that, for every $X \in L(\HH),$ $\vert (KX-B_0') \vert^2 = (KX-B_0')^*(KX-B_0') = \vert TX \vert^2 + \vert VX-B_0 \vert^2.$
	Then, for every $X \in L(\HH)$
	$$\Vert TX_0 \Vert_{2}^2 + \Vert VX_0-B_0 \Vert_2^2 = \tr (\vert (KX_0-B_0') \vert^2)  \leq \tr (\vert (KX-B_0') \vert^2) =\Vert TX \Vert_{2}^2 + \Vert VX-B_0 \Vert_2^2,$$ see Proposition \ref{Prop Nashed}.
	Thus, \eqref{dos} admits a solution.
	
	To prove the converse, consider  
	$F_2: L(\HH) \rightarrow \mathbb{R}^{+},$ 
	$$\begin{array}{lll}
	F_2(X)&=&\Vert TX \Vert_{2}^2 + \Vert VX-B_0 \Vert_2^2= \Vert T X \Vert_{2}^2 + \tr ((VX-B_0)^*(VX-B_0))\\
	&=&  \Vert  T                X\Vert_{2}^2+  \Vert V X\Vert_{2}^2-2Re \ [\tr(B_0^*VX)]+ \Vert  B_0\Vert_{2}^2           \\
	&=&\Vert \vert T\vert X 
	\Vert_{2}^2+\Vert \vert V\vert X\Vert_{2}^2-2 Re \  [\tr(B_0^*VX)]+ \Vert  B_0\Vert_{2}^2\\
	\end{array}
	$$
	By Theorem \ref{TeoD}, $F_2$ has a Fr\'echet derivative and, furthermore, for every $X, Y \in L(\HH)$ 
	$$\begin{array}{lll}
	DF_2 (X)(Y) &=& DG_2 (\vert T\vert X)(\vert T\vert Y) +DG_2 (\vert V \vert X)(\vert V \vert Y)-2Re \ [\tr(B_0^*VY)]\\
	&=& 	2 Re \ [ \tr ( (\vert T \vert X_0)^*\vert T \vert Y) ]+ 2 Re \ [ \tr( (\vert V \vert X_0)^* \vert V \vert Y)] -2 Re \ [ \tr (B_0^*VY)],
	\end{array}
	$$
	where $G_2(X)= \Vert X \Vert_2^2=\tr (X^*X).$
	
	Suppose that $X_0 \in L(\HH)$ is a global minimum  of $\Vert TX \Vert_{2}^2 + \Vert VX-B_0 \Vert_2^2.$ Then $X_0$ is a global minimum of $F_2$ and, since $F_2$ is a $\mc{C}^1$-function
	$$DF_2(X_0)(Y)= 0, \mbox{ for every } Y \in L(\HH).$$
	Then, for every $Y \in L(\HH),$ 
	$$Re \ [ \tr (( (\vert T \vert X_0)^*\vert T \vert  + (\vert V \vert X_0)^* \vert V \vert  - B_0^*V) \ Y)]=0.$$
	Then,  it follows that
	$$X_0^*\vert T\vert^2 +X_0^*\vert V \vert ^2-B_0^*V=0$$ or, equivalently  
	$$(T^*T+V^*V)X_0=V^{*}B_0.$$ 
\end{dem}

To study the existence of solutions of inconsistent linear systems under seminorms defined by positive semidefinite matrices, Mitra defined the optimal inverses for matrices \cite{[Mit]}. In \cite{[CorFonMae16]}, Mitra's concept was extended to Hilbert spaces:

\begin{Def}
	Given operators $A \in L(\HH,\mc{F})$ and $W\in L(\mc{F}\oplus \HH)^+,$ a \emph{$W$-optimal inverse} of $A$ is an operator $G \in L(\mc{F},\HH)$ such that 
	$$\Vert \left(
	\begin{array}{cc}
	AGf-f \\
	Gf \\
	\end{array}
	\right) \Vert_{W} = \underset{h \in \HH}{\min} \ \Vert \left( \begin{array}{cc} 
	Ah-f \\
	h \\
	\end{array}
	\right)\Vert_{W} ,$$ for every $f \in \mc{F}.$ Here $\Vert \cdot \Vert_{W}$ denotes the seminorm defined by $W:$ $\Vert \left( \begin{array}{cc} 
	f \\
	h \\
	\end{array}
	\right) \Vert_W=\Vert W^{1/2}  \left( \begin{array}{cc} 
	f \\
	h \\
	\end{array}
	\right) \Vert.$
\end{Def}


\bigskip
Consider  $W$ with the following block form 
\begin{equation} \label{AblockformB}
\left( \begin{array}{cc} 
W_{11} & W_{12} \\
W_{12}^* & W_{22}\\
\end{array}
\right),
\end{equation}
where $W_{11} \in L(\mc{F})^+,$ $W_{22} \in L(\HH)^+$ and $W_{12} \in L(\HH,\mc{F}).$ By \cite[Theorem 2.1]{[CorFonMae16]} and \cite[Theorem 4.2]{[Mit]}, $A \in L(\HH,\mc{F})$ admits a $W$-optimal inverse if and only if the equation
\begin{equation} \label{AoptinvB}
(A^*W_{11}A+A^*W_{12}+W_{12}^*A+W_{22})X=A^*W_{11}+W_{12}^*
\end{equation}
admits a solution. In this case, the set of $W$-optimal inverses of $A$ is the set of solutions of \eqref{AoptinvB}.

\vspace{0,3cm}
The following result relates the existence of a $W$-optimal inverse to the existence of a solution of  \eqref{dos}. Some equivalences of the next proposition were proven in \cite[Theorem 4.2]{[CorFonMae16]} for $V \in L(\HH,\mc{F})$ with closed range. The proofs of such equivalences are included in order to remark that the range of $V$ need not be closed.

\begin{thm} \label{cor31} Let $T \in S_2(\HH,\mc{E})$ and  $V \in S_2(\HH, \mc{F}).$ Then the following are equivalent:
	\begin{itemize}
		\item [i)] there exists $\underset{X \in L(\HH)}{\min} \ (\Vert TX \Vert_{2}^2 + \Vert VX-B_0 \Vert_2^2)$ for every $B_0 \in S_2(\HH, \mc{F}),$ 
		\item [ii)] $R(V^*) \subseteq R(T^*T+V^*V),$
		\item [iii)] there exists $\underset{h \in \HH}{\min} \ (\| Th\|^2+\|Vh-f_0\|^2)$ for every $f_0\in \mathcal{F},$
		\item[iv)] $V$ admits a $\left( \begin{array}{cc} 
		I & 0 \\
		0 & T^*T\\
		\end{array}
		\right)$-optimal inverse,
		\item [v)] there exists a bounded global solution of the classical smoothing problem \eqref{CSP}.
	\end{itemize}		
	If $R(V)$ is closed, conditions $i)$ to $v)$ are also equivalent to
	\begin{itemize}
		\item [vi)] the pair $(T^*T,N(V))$ is compatible.
	\end{itemize}
\end{thm}
\begin{dem} 
	$i) \Leftrightarrow ii):$ Suppose that \eqref{dos} has a minimum for every $B_0 \in S_2(\HH, \mc{F}).$ Then, by Proposition \ref{teo3} and Douglas' Lemma, $R(V^*B_0) \subseteq R(T^*T+V^*V).$ 
	Consider $f_0 \in \mc{F}$, then  there exists $B_0 \in S_2(\HH, \mc{F})$ such that $f_0=B_0x,$ for some $x\in \mathcal H$. Therefore
	$$
	V^*f_0=V^*B_0x \in R(V^*B_0)\subseteq R(T^*T+V^*V).
	$$
	Hence $R(V^*) \subseteq R(T^*T+V^*V).$
	Conversely, suppose that  $R(V^*) \subseteq R(T^*T+V^*V)$ and consider  $B_0 \in S_2(\HH, \mc{F}).$ Then  $R(V^*B_0)\subseteq R(V^*) \subseteq R(T^*T+V^*V)$.  Hence, by Douglas' Lemma and Proposition \ref{teo3}, \eqref{dos} has a solution for every $B_0 \in S_2(\HH, \mc{F}).$ 
	
	$ii) \Leftrightarrow iii):$ Let $K$ be as in \eqref{eqK}. Given $f_0\in \mc{F}$, it holds that $
	\underset{h \in \HH}{\min} \ (\| Th\|^2+\|Vh-f_0\|^2)=\underset{h \in \HH}{\min} \ \|Kh-(0,f_0)\|^2
	$ exists if and only if the normal equation $K^*K h=K^*(0, f_0)$ has a solution; equivalently $V^*f_0=(T^*T+V^*V)h$ has a solution. Therefore, $
	\underset{h \in \HH}{\min}\ (\| Th\|^2+\|Vh-f_0\|^2)
	$ exists for every $f_0\in \mathcal{H}$ if and only if $R(V^*)\subseteq R(T^*T+V^*V)$.
	
	$ii) \Leftrightarrow iv):$ It follows by \eqref{AoptinvB} and Douglas' Lemma.
	
	$iv) \Leftrightarrow v):$ Consider the inner product and the associated norm on $\mc{F} \oplus \HH$ as in \eqref{pi}. Then $G \in L(\mc{F},\HH)$ is a $\left( \begin{array}{cc} 
	I & 0 \\
	0 & T^*T\\
	\end{array}
	\right)$-optimal inverse of $V$ if and only if for every $f_0 \in \mathcal F$,  $\|VGf_0-f_0\|^2+\|Gf_0\|_{T^*T}^2\leq  \|Vh-f_0\|^2+\|h\|_{T^*T}^2 $ for every $h \in \mathcal H$ or, equivalently, $\|VGf_0-f_0\|^2+\|TGf_0\|^2\leq  \|Vh-f_0\|^2+\|Th\|^2 $ for every $h\in \mathcal H$, that is, $G$ is a bounded global solution of \eqref{CSP}.
	
	$ii) \Leftrightarrow vi):$ Suppose that $R(V)$ is closed, then $R(V^*V)=R(V^*)$. It can be seen that  $R(V^*)=R(V^*V)\subseteq R(T^*T+V^*V)$ if and only if $R(T^*T+V^*V)=R(T^*T)+R(V^*V)$. But this is equivalent to  $(T^*T, N(V))$ being compatible, see \cite[Theorem 3.2]{[CorFonMae16]}.
\end{dem}

\bigskip
In Theorem \ref{cor31}, it was proved that the existence of a bounded global solution of the classical smoothing problem \eqref{CSP} is equivalent to the existence of a  $W$-optimal inverse for the weight $\left( \begin{array}{cc} 
I & 0 \\
0 & T^*T\\
\end{array}
\right).$ In a similar way, in Theorem \ref{sol global sii  A-inversa}, the equivalence between the existence of bounded global solutions for \eqref{WLSP} and the  existence of $W$-inverses was stated for a positive weight $W$. Motivated by this relation, in what follows we are interested in comparing $W$-inverses to $W$-optimal inverses. We begin with the following lemma.

\begin{lema} \label{thmAopAinv} Let $W\in L(\mc{F}\oplus \HH)^+$ with block form as in \eqref{AblockformB}, $A \in L(\HH,\mc{F})$ and $\hat{A}\in L(\HH,\mc{F}\oplus \HH)$ be defined by $\hat{A}h=(Ah,h)$. Then, there exists a $W$-inverse of $\hat{A}$ if and only if there exists a $W$-optimal inverse of $A$ and the equation 
	\begin{equation} \label{Aopt2}
	(A^*W_{11}A+A^*W_{12}+W_{12}^*A+W_{22})X=A^*W_{12}+W_{22}
	\end{equation}
	admits a solution.
\end{lema}

\begin{dem} Suppose that $Z \in L(\mc{F}\oplus \HH,\HH)$ is a $W$-inverse of $\hat{A},$ then
	$\hat{A}^*W\hat{A}Z(f,h)=\hat{A}^*W(f,h)$ for every $(f,h) \in \mc{F} \oplus \HH,$ i.e.,
	$$(A^*W_{11}A+A^*W_{12}+W_{12}^*A+W_{22})Z(f,h)=A^*W_{11}f+A^*W_{12}h+W_{12}^*f+W_{22}h,$$ 
	for every $(f,h) \in  \mc{F} \oplus \HH.$
	In particular, if $h=0,$ then 
	$$(A^*W_{11}A+A^*W_{12}+W_{12}^*A+W_{22})Z(f,0)=A^*W_{11}f+W_{12}^*f, \mbox{ for every } f \in \mc{F}.$$
	Therefore $Z_1(f) :=Z(f,0)$ is a $W$-optimal inverse of $A$.
	In the same way, if $f=0,$ then 
	$$(A^*W_{11}A+A^*W_{12}+A_{12}^*A+W_{22})Z(0,h)=A^*W_{12}h+W_{22}h, \mbox{ for every } h \in \HH.$$
	Therefore $Z_2(h) :=Z(0,h)$ is a solution of \eqref{Aopt2}.
	
	Conversely, suppose that $Z_1 \in L(\mc{F},\HH)$ is a $W$-optimal inverse of $A$ and $Z_2 \in L(\HH)$ is a solution of \eqref{Aopt2}. Let $Z: \mc{F} \oplus \HH \ra \HH$ be defined by $Z(f,h):=Z_1(f)+Z_2(h).$ Then, clearly $$\hat{A}^*W\hat{A}Z=\hat{A}^*W.$$
	Also, 
	\begin{align*}
	\Vert Z(f,h) \Vert^2 &= \Vert Z_1(f) + Z_2(h) \Vert^2 \leq (\Vert Z_1(f) \Vert + \Vert Z_2(h) \Vert)^2  \leq (\Vert Z_1 \Vert \Vert f \Vert + \Vert Z_2 \Vert \Vert h \Vert)^2\\
	&\leq(\max\{\Vert Z_1 \Vert, \Vert Z_2 \Vert\})^2 (\Vert f \Vert + \Vert h \Vert)^2 \leq 2 (\max\{\Vert Z_1 \Vert, \Vert Z_2 \Vert\})^2 (\Vert f \Vert^2 + \Vert h \Vert^2)\\
	&=2 (\max\{\Vert Z_1 \Vert, \Vert Z_2 \Vert\})^2 \Vert (f,h)\Vert^2.
	\end{align*}
 Therefore $Z \in L(\mc{F}\oplus \HH, \HH)$ and $Z$ is a $W$-inverse of $\hat{A}.$
\end{dem}

The next proposition shows that certain optimal inverses can be seen as the weighted inverse of an associated operator. 
\begin{prop} Let $A \in L(\HH,\mc{F}),$ $w \in L(\HH)^+$ and $\hat{A}\in L(\HH,\mc{F}\oplus \HH)$ be defined by $\hat{A}h=(Ah,h)$. Then, there exists a $\left( \begin{array}{cc} 
	I & 0 \\
	0 &w \\
	\end{array}
	\right)$-inverse of $\hat{A}$ if and only if there exists a $\left( \begin{array}{cc} 
	I & 0 \\
	0 & w\\
	\end{array}
	\right)$-optimal inverse of $A.$ 
\end{prop}

\begin{dem}
	Suppose that there exists a $\left( \begin{array}{cc} 
	I & 0 \\
	0 & w\\
	\end{array}
	\right)$-optimal inverse of $A.$ By \eqref{AoptinvB} and Douglas' Lemma, $R(A^*) \subseteq R(A^*A+w).$ This is $R(A^*A)\subseteq R(A^*A+w)$ and, this is equivalent to $R(A^*A)+R(w) = R(A^*A+w).$ Hence, $R(w) \subseteq R(A^*A+w)$ and by Douglas's Lemma, the equation $(A^*A+w)X=w$ admits a solution. Therefore, by Lemma \ref{thmAopAinv}, $\hat{A}$ admits a  $\left( \begin{array}{cc} 
	I & 0 \\
	0 & w\\
	\end{array}
	\right)$-inverse. The converse follows by Lemma \ref{thmAopAinv}.
\end{dem}

\section*{Acknowledgements}
We thank the anonymous referees for carefully reading our manuscript and helping us to improve
this paper with several useful comments.

Maximiliano Contino was supported by CONICET PIP 0168. Guillermina Fongi was partially supported by ANPCyT PICT 2017-0883.

\end{document}